\documentclass[10pt]{article}
\usepackage{amsmath,amssymb,amsthm,amscd}
\numberwithin{equation}{section}

\newtheorem{prop}{Proposition}[section]
\newtheorem{theorem}[prop]{Theorem}

\newtheorem{corollary}[prop]{Corollary}
\newtheorem{remark}[prop]{Remark}

\newtheorem{acknowledge}[prop]{Acknowledgment}

\def\<{\langle}
\def\>{\rangle}
\def\({\left(}
\def\){\right)}

\def\p{\partial}
\def\Ric{{\rm Ric}}

\begin{document}

\title{K\"ahler-Ricci Flow with Degenerate
Initial Class}
\author{Zhou Zhang \footnote{Research supported
in part by National Science Foundation grants
DMS-0904760.} \\
University of Michigan \\
at Ann Arbor}
\date{}
\maketitle
\begin{abstract}

\noindent
In \cite{chen-tian-z}, the weak K\"ahler-Ricci
flow was introduced for various geometric
motivations. In the current work, we take
further consideration on setting up the weak
flow. Namely, the initial class is allowed to
be no longer K\"ahler. The convergence as
$t\to 0^+$ is of great importance to study
for this topic.

\end{abstract}

\section{Motivation and Set-up}

K\"ahler-Ricci flow, the complex version
of Ricci flow, has been under intensive
study over the past twenty some years.
In \cite{tian02} and more recently \cite
{tian-survey}, G. Tian proposed the
intriguing program of constructing globally
existing (weak) K\"ahler-Ricci flow with
canonical (singular) limit at infinity and
applying it to the study of general
algebraic manifold.

Generally speaking, one should expect
the classic K\"ahler-Ricci flow to encounter
singularity at some finite time which is
completely decided by cohomology
information according to the optimal
existence in \cite{t-znote}. Just as what
people wanted to do and have had
successes in some cases for Ricci flow,
surgery on the underlying manifold
should be expected. For K\"ahler-Ricci
flow, we naturally want the surgery to
have more flavor in algebraic geometry.
For surface of general type, we only
need the blow-down of $(-1)$-curves
to apply the construction in \cite
{chen-tian-z} to push the flow through
finite time singularities. The degenerate
class at the singularity time would
become K\"ahler for the new manifold
because the $(-1)$-curves causing the
cohomology degeneration have been
crushed to points. Things can get
significantly more complicated for higher
dimensional manifold. In (complex)
dimension $3$, flips are involved.
Simply speaking, one needs to blow
up the manifold and then blow down.
Naturally, we could expect the
transformation of the degenerate class
is still not K\"ahler. In this note, we want
to say this is not a problem if formally
the K\"ahler-Ricci flow is instantly taking
the class into the K"ahler cone of the
new manifold. As in \cite{chen-tian-z},
short time existence is the topic. In the
following, the precise problem under
consideration is stated and we set up
the a priori weak flow following the
same idea as in the previous work.

\vspace{0.1in}

Let $X$ be a closed K\"ahler manifold with
$dim_{\mathbb{C}}X=n\geqslant 2$. We
consider the following version of K\"ahler-Ricci
flow over $X$,
\begin{equation}
\label{eq:krf}
\frac{\partial\widetilde\omega_t}{\partial t}
=-{\rm Ric}(\widetilde\omega_t)-\widetilde
\omega_t, ~~~~\widetilde\omega_0=\omega
_0.
\end{equation}
More importantly, the weak initial data
$$\omega_0=\omega+\sqrt{-1}\p\bar\p
v$$
where $\omega$ is a real, smooth, and
closed $(1, 1)$-form with $[\omega]$
being nef. (i.e. numerically effective, in
other words, on the boundary of the
K\"ahler cone of $X$), and $v\in PSH_
\omega(X)\cap L^\infty(X)$ where $v\in
PSH_\omega(X)$ means $v+\sqrt{-1}\p
\bar\p v\geqslant 0$ weakly (i.e. in the
sense of distribution).

\begin{remark}

The choice of K\"ahler-Ricci in this version
is not essential. Our discussion is even valid
for other unconventional K\"ahler-Ricci type
of flows as in \cite{t-znote}.

\end{remark}

The main additional assumption for this work
is the following.

{\bf Suppose $[\omega_0]=[\omega]$ is on
the boundary of the K\"ahler cone of $X$
(in the cohomology space $H^{1, 1}(X,
\mathbb{C})\cap H^2(X, \mathbb{R})$).
The ray starting from $[\omega]$ and in
the direction towards the canonical class
of $X$, $K_X$ enters the K\"ahler cone
instantly.}

Clearly, there is no need for $K_X$ to be
K\"ahler for this to happen.

There are other motivations to study this
case besides defining weak flow to realize
Tian's program as mentioned before. In
general, there is this understanding that
the existence of K\"ahler-Ricci flow to be
decided completely by the cohomology
information from formal ODE consideration.
The main theorem below would strengthen
this philosophic point of view.

\begin{theorem}

\label{th:main}
In the above setting, one can define a
unique weak K\"ahler-Ricci flow from
the approximation construction which
becomes smooth instantly and satisfying
(\ref{eq:krf}).

\end{theorem}

Formally, we see $[\widetilde\omega_t]=
[\omega_t]\in H^{1, 1}(X, \mathbb{C})\cap
H^2(X, \mathbb{R})$ where
$$\omega_t=\omega_\infty+e^{-t}(\omega-
\omega_\infty)$$
with $\omega_\infty=-\Ric(\Omega):=-\sqrt
{-1}\p\bar\p\log\(\frac{\Omega}{(\sqrt{-1})^n
dz_1\wedge d\bar z_1\wedge\cdots\wedge
dz_n\wedge d\bar z_n}\)$ in a local
coordinate system $\{z_1, \cdots, z_n\}$.
The main additional assumption above simply
means
$$[\omega_t]=e^{-t}[\omega]+(1-e^{-t})K_X$$
is K\"ahler for $t\in (0, T)$ for some $T>0$.
It is now routine to see at least formally (\ref
{eq:krf}) would be equivalent to the following
evolution equation for a space-time function
$u$,
\begin{equation}
\label{eq:skrf}
\frac{\p u}{\p t}=\log\frac{(\omega_t+\sqrt{-1}
\p\bar\p u)^n}{\Omega}-u, ~~~~u(\cdot, 0)=v
\end{equation}
with the understanding of
$\widetilde\omega_t=\omega_t+\sqrt{-1}\p
\bar\p u$.

Just as for classic smooth K\"ahler-Ricci flow,
we'll focus on defining weak version of (\ref
{eq:skrf}) instead of (\ref{eq:krf}). Their
equivalence in the category of smooth
objects would make the weak version for
(\ref{eq:skrf}) naturally the weak version for
(\ref{eq:krf}) in sight of the smoothing effect
in Theorem \ref{th:main}.

\vspace{0.1in}

It's time to describe the {\bf approximation
construction} mentioned in Theorem \ref
{th:main}, which is similar to what has been
applied in \cite{chen-tian-z} except that now
$[\omega]$ is no longer K\"ahler. The idea
is to find approximation of the initial data,
use them as initial data to get a sequence
of flows and finally take limit of the flows.
The detail is as follows.

Take some K\"ahler metric $\omega_1$ over $X$.
For any $\epsilon\geqslant 0$, set $\omega(\epsilon)
=\omega+\epsilon\omega_1$ and $\omega_t
(\epsilon)=\omega_\infty+e^{-t}\(\omega(\epsilon)-
\omega_\infty\)$. Using the regularization result in
\cite{blo-koj}, one has for any sequence $\{\epsilon
_j\}$ decreasing to $0$ as $j\to\infty$, a sequence
of functions $\{v_j\}$ with $v_j\in C^\infty(X)$ and
$\omega(\epsilon_j)+\sqrt{-1}\p\bar\p v_j>0$,
decreasing to $v$ accordingly. Then one considers
the K\"ahler-Ricci flows,
\begin{equation}
\label{eq:approxi-krf}
\frac{\partial\widetilde\omega_t(\epsilon_j)}{\partial
t}=-{\rm Ric}\(\widetilde\omega_t(\epsilon_j)\)-
\widetilde\omega_t(\epsilon_j), ~~~~\widetilde
\omega_0(\epsilon)=\omega(\epsilon_j)+\sqrt{-1}
\p\bar\p v_j.
\end{equation}
At the level of metric potential, we have
\begin{equation}
\label{eq:approxi-skrf}
\frac{\p u_j}{\p t}=\log\frac{(\omega_t(\epsilon_j)
+\sqrt{-1}\p\bar\p u_j)^n}{\Omega}-u_j, ~~~~u_j
(\cdot, 0)=v_j.
\end{equation}
They are in the classic setting of K\"ahler-Ricci
flow. By choosing the $T$ (before (\ref{eq:skrf}))
properly, all these flows for $j\gg 1$ (i.e.
$\epsilon_j$ sufficiently small) would exists for
$t\in [0, T)$ from cohomology consideration
\footnote{As in \cite{chen-tian-z}, it's the short
time existence that we are interested in, i.e.
the small interval near $t=0$.}.

In sight of $v_j$ and $\omega_t(\epsilon_j)$
decreasing to $v$ and $\omega_t$ as $j\to
\infty$. Applying Maximum Principle, one can
see $u_j$ is also decreasing as $j\to\infty$. In
principle, this would allow us to take and get
a limit for each $t\in [0, T)$, $u(\cdot, t)\in
PSH_{\omega_t}(X)$, which is the weak
flow wanted. For this to be true literally, one
needs to make sure for each such $t$, the
decreasing limit of $u_j(\cdot, t)$ won't be
$-\infty$ uniformly over $X$. At the initial
time, this is obviously the case. For $t\in
(0, T)$, this would be the case as seen in
later by applying Kolodziej's $L^\infty
$-estimate (as in \cite{kojnotes}) \footnote
{One can also achieve this using the classic
PDE argument involving Moser Iteration.}.

Of course, one needs to make sure the limit,
which is at this moment just a family of
positive $(1,1)$-current with parameter $t$,
is more classic. This clearly boils down to
obtain uniform estimates for $u_j$.

\vspace{0.1in}

Before heading into the search for
those uniform estimates, let's justify
the uniqueness statement of Theorem
\ref{th:main}. This is again a fairly
direct application of Maximum Principle.

To begin with, let's observe that the
decreasing convergence of $v_j\to
v$ can be strengthened to strictly
decreasing convergence without
changing the limit. Suppose $\{v_j
\}$ is only a decreasing sequence,
then $\{v_j+\frac{1}{j}\}$ is strictly
decreasing with the same limit.
Clearly, the affect on the solution
of (\ref{eq:approxi-skrf}), $u_j$ is
negligible as $j\to\infty$. Also, the
decreasing limit won't be affected
by taking subsequence.

Secondly, we see the choice of sequence
$\{\epsilon_j\}$ won't affect the limit (i.e.
the weak flow). Let's take two strictly
decreasing sequences, $\{v_j\}$ and $\{
v_\alpha\}$ in the construction before
(\ref{eq:approxi-krf}). Since for each $j$,
$v_j>v$ and $v_\alpha$ decreases to
$v$, by Dini's Theorem, $v_\alpha<v_j$
for $\alpha$ sufficiently large. The other
direction is also right. So by taking
subsequences, still denoted by $\{v_{j_a}
\}$ and $\{v_{\alpha_b}\}$, we have
$$\cdots<v_{\alpha_2}<v_{j_2}<v_
{\alpha_1}<v_{j_1}.$$
Applying Maximum Principle to (\ref
{eq:approxi-skrf}), one has
$$\cdots<u_{\alpha_2}<u_{j_2}<u_
{\alpha_1}<u_{j_1},$$
and so they have the same limit.

Now we take care of the general case. In
the construction before (\ref{eq:approxi-krf}),
suppose we have chosen different $\omega
_1$'s with different strictly decreasing
sequences $\{v_j\}$. Then we can use
the argument in the above situation and
also make sure the corresponding $\omega
_t(\epsilon_j)$'s have the same comparison
relation. Thus Maximum Principle would
still give the same kind of comparison for
solutions of (\ref{eq:approxi-skrf}). Hence
the limit would still be the same.

\begin{remark}

It's not hard to make this construction even
more flexible, and so far, this is pretty much
the only way to come up with a reasonable
weak flow. So even without any description
of the situation as $t\to 0^+$, it is not so
artificial to call this a weak flow initiating
from $\omega_0=\omega+\sqrt{-1}\p\bar
\p v$.

\end{remark}

\section{Local General Estimates}

Now we begin the search for estimates uniform
for all approximation flows (i.e. $\epsilon>0$
where we have such a flow). For simplicity, we'll
{\bf omit $j$ and $\epsilon_j$ in the notations}
below which would unfortunately make (\ref
{eq:approxi-krf}) and (\ref{eq:approxi-skrf}) look
exactly like (\ref{eq:krf}) and (\ref{eq:skrf})
respectively. However this would help us
keeping in mind about the degeneracy of the
background form, and so is not such a terrible
choice considering what we are trying to do.

{\bf Note:} $C$ below would stand for fixed
positive constant which might be different
at places. Its dependence on other constants
should be clear from the context.

\vspace{0.1in}

Clearly, $u\leqslant C$ by Maximum Principle
for (\ref{eq:approxi-skrf}). This is for all time.
For the other estimates, the idea is trying to
eliminate the affect of the initial data as
completely as possible because we don't
have control of the initial data except for $L^
\infty$-bound of the potential. Also recall that
we only need estimates for short time.

\vspace{0.1in}

Notice that in (\ref{eq:approxi-skrf}), the initial
value of the background form, $\omega$ might
not be non-negative. One can actually make
better use of that $\omega+\sqrt{-1}\p\bar\p
v>0$ by looking at the following evolution at
the level of metric potential for the same flow
(\ref{eq:approxi-krf}),
\begin{equation}
\label{eq:approxi-skrf2}
\frac{\p \phi}{\p t}=\log\frac{(\widehat\omega_t
+\sqrt{-1}\p\bar\p \phi)^n}{\Omega}-\phi, ~~~~
\phi(\cdot, 0)=0,
\end{equation}
where $\widehat\omega_t=\omega_\infty+
e^{-t}(\omega+\sqrt{-1}\p\bar\p v-\omega_
\infty)$. It's easy to see the relation between
the solutions of (\ref{eq:approxi-skrf}) and
(\ref{eq:approxi-skrf2}) is $u=\phi+e^t\cdot
v$. So $\phi\leqslant C$, which is not clear
by applying Maximum Principle to (\ref
{eq:approxi-skrf2}) because of the lack of
uniform control for $\omega+\sqrt{-1}\p\bar
\p v$ as form.

This following equation is obtained by taking
$t$-derivative of (\ref{eq:approxi-skrf2}) and
making some transformations
\begin{equation}
\frac{\partial}{\partial t}\((e^t-1)\frac{\partial
\phi}{\partial t}-\phi\)=\Delta_{\widetilde
\omega_t}\((e^t-1)\frac{\partial\phi}{\partial
t}-\phi\)+n-\<\widetilde\omega_t,\omega+
\sqrt{-1}\p\bar\p v\>. \nonumber
\end{equation}
Since $\omega+\sqrt{-1}\p\bar\p v>0$, applying
Maximum Principle and noticing the lower
bound of the initial value and the uniform
upper bound of $\phi$, we have for $t\in (0, T)$,
$$\frac{\partial\phi}{\partial t}\leqslant\frac{C}
{e^t-1},$$
which gives the following bound of $\frac
{\p u}{\p t}$ since $\frac{\p u}{\p t}=\frac{\p\phi}
{\p t}+e^t\cdot v$,
$$\frac{\partial u}{\partial t}\leqslant\frac{C}
{e^t-1}.$$

For any $t\in (0, T)$, we have the background
form $[\omega_t(\epsilon)]$ being uniformly
K\"ahler, i.e. the small interval corresponding
to $\epsilon$ is in the K\"ahler cone. Together
with the above upper bound for $\frac{\p u}{\p
t}$, one can apply Ko\l odziej's $L^\infty
$-estimate (as in \cite{kojnotes}) for (\ref
{eq:approxi-skrf}) in the form of
$$(\omega_t+\sqrt{-1}\p\bar\p u)^n=e^{\frac
{\p u}{\p t}+u}\Omega$$
to achieve the $L^\infty$-bound for $u$. So
now we also have $u(\cdot, t)\geqslant -C(t)$
with $C(t)$ finite for $t\in (0, T)$. In fact, we
know by the result in \cite{koj06} that $u(\cdot,
t)$ is H\"older continuous for these $t$'s
\footnote{The H\"older exponent will also
depend on $t$.}.

\begin{remark}

\label{re:about-measure}
The original results on $L^\infty$-estimate (as
in \cite{kojnotes}, \cite{t-znote} and \cite{zzo})
are usually stated for Monge-Amp\`ere
equation in the form $(\omega+\sqrt{-1}\partial
\bar\partial u)^n=f\cdot\Omega$ where
$[\omega]$ might be degenerate, $f\geqslant
0$ is in some $L^{p>1}$-space and $\Omega$
is a (non-degenerate) smooth volume form.
There are more than one way to translate this
when applying to the equation with $e^u$ on
the right hand side.

Method I: get $L^p$ bound for the measure
$f\cdot e^u\Omega$, then one knows the
normalized $u$ would be bounded from
the original result. In other words, $u$ only
takes value in some interval with well
controlled length. Then the upper bound of
$u$, which usually comes from direct
Maximum Principle argument, together with
the upper bound for $f\cdot\Omega$ which
guarantees $u$ can not take too small value
all over in sight of the total volume having
a lower bound, would provide the bound for
$u$ itself.

Method II: get $L^p$ bound for the measure
$f\cdot\Omega$, then consider the equation
$(\omega+\sqrt{-1}\partial\bar\partial w)^n=
Cf\cdot\Omega$. The idea is to apply
Maximum Principle to the quotient of these
two equations. In order to control the
(normalized) solution $w$ for this auxiliary
equation, one needs to control the constant
$C$ (from above), which means we need a
lower bound for the total volume for the
measure $f\cdot\Omega$. Again this can
be achieved from the upper bound for $u$.

It's not hard to see that these two methods
are merely different combinations of the
same set of information.

\end{remark}

Since for any $t\in [\lambda_1, \lambda_2]
\subset (0, T)$, $[\omega_t]$ is uniformly
K\"ahler (for any approximation flow), by
properly choosing $\Omega$ (and so
$\omega_\infty$) \footnote{This is only
going to cause difference for the evolution
equations at the level of metric potential
similar to that between (\ref{eq:approxi-skrf})
and (\ref{eq:approxi-skrf2}).}, one has $
\omega_t$ being uniform as K\"ahler metric
for $t\in [\lambda_1, \lambda_2]$ \footnote
{In fact, one only needs to take care of the
case $\epsilon=0$ to achieve this.}. Clearly,
$(0, T)$ can be exhausted by such closed
intervals. Of course, we only care for the
end towards $t=0$.

Now let's translate the time to make
$\lambda_1$ the new initial time. From
the discussion before, we have uniform
bounds from both sides for $u$ and the
uniform upper bound for $\frac{\p u}{\p
t}$. By taking $t$-derivative for (\ref
{eq:approxi-skrf}) and making some
transformations, we have the two
equations below
\begin{equation}
\label{eq:1}
\frac{\partial }{\partial t }\((e^t-1)\frac
{\partial u }{\partial t }\)=\Delta_{\widetilde
{\omega}_t }\((e^t-1)\frac{\partial u }
{\partial t }\)-(1-e^{-t})\<\widetilde
{\omega }_t, \omega-\omega_\infty
\>+\frac{\partial u }{\partial t },
\end{equation}
\begin{equation}
\label{eq:2}
\frac{\partial}{\partial t}\((e^t-1)\frac
{\partial u}{\partial t}-u\)=\Delta_{
\widetilde\omega_t}\((e^t-1)\frac{
\partial u}{\partial t}-u\)+n-\<\widetilde
\omega_t,\omega\>.
\end{equation}
Notice that the $\omega$ here is indeed
$\omega_{\lambda_1}$ and is uniform
as K\"ahler metric for all approximation
flows.

In the small time interval $[0, \lambda_2-
\lambda_1]$ (after translation), we have
made sure that $\omega_t>0$, which is
$$\omega_\infty+e^{-t}(\omega-\omega_
\infty)=\omega-(1-e^{-t})(\omega-\omega
_\infty)>0,$$
so one can choose $\delta\in (0, 1)$ such
that for these $t$'s,
$$\delta\omega-(1-e^{-t})(\omega-\omega
_\infty)>0.$$

Use this $\delta$ to multiply (\ref{eq:2})
and take difference with (\ref{eq:1}) to
arrive at
\begin{equation}
\begin{split}
\frac{\partial }{\partial t }\((1-\delta)(e^t-
1)\frac{\partial u }{\partial t}+\delta u\)
&= \Delta_{\widetilde{\omega}_t }\((1-
\delta)(e^t-1)\frac{\partial u }{\partial t}
+\delta u\)\\
&~~~~ +\<\widetilde{\omega }_t, \delta
\omega-(1-e^{-t})(\omega-\omega_\infty)
\>+\frac{\partial u }{\partial t }-n\delta.
\nonumber
\end{split}
\end{equation}

Consider the minimum value point of
the term $(1-\delta)(e^t-1)\frac{\partial
u}{\partial t}+\delta u$. If it is not at the
(new) initial time, then at that point, we
have
\begin{equation}
\begin{split}
\<\widetilde{\omega }_t, \delta\omega
-(1-e^{-t})(\omega-\omega_\infty)\>
&\geqslant n\cdot\(\frac{\(\delta\omega
-(1-e^{-t})(\omega-\omega_\infty)\)^n}
{\widetilde\omega^n_t}\)^{\frac{1}{n}} \\
&= n\cdot\(\frac{\(\delta\omega-(1-e^{-t})
(\omega-\omega_\infty)\)^n}{e^{\frac
{\p u}{\p t}+u}\Omega}\)^{\frac{1}{n}} \\
&\geqslant Ce^{-\frac{1}{n}\frac{\p u}
{\p t}}
\end{split}
\end{equation}
where $u\leqslant C$ is applied in the
last step, and so one conclude
$$C\geqslant Ce^{-\frac{\partial u}
{\partial t}}+\frac{\partial u}{\partial t},$$
which gives $\frac{\partial u}{\partial t}
\geqslant -C$ at that point. So $(1-\delta)
(e^t-1)\frac{\partial u}{\partial t}+\delta u
\geqslant -C$ at that point in sight of the
lower bound of $u$ from previous
argument. We conclude that
$$(1-\delta)(e^t-1)\frac{\partial u}{\partial
t}+\delta u\geqslant -C$$
for the space-time, and so
$$\frac{\partial u }{\partial t }\geqslant
-\frac{C}{e^t-1}.$$
Remember that the time has been
translated.

\vspace{0.1in}

So far, we have obtained the $L^\infty$-bounds
for both $u$ and $\frac{\partial u}{\partial t}$
locally away from the initial time. Only the upper
bound of $u$ is uniform for all time.

The second and higher order estimates can
be carried through as in Subsection 3.2 of
\cite{chen-tian-z} because the translation of
time would make the background form
K\"ahler. Hence we conclude that the weak
flow defined in Section 1 becomes smooth
instantly. The proof of Theorem \ref{th:main}
is thus finished.

\begin{remark}

The situation as $t\to 0^+$, which is indeed
the only "weak" spot of the flow, needs further
consideration just as in \cite{chen-tian-z}.
Since most estimates achieved up to this
point are only local away from the initial
time, the control of the situation near $0$
at this moment is very weak. In fact, strange
things can happen. For example, $[\omega_
0]$ might have $0$ volume (being collapsed),
but the volume would becomes positive
instantly.

\end{remark}

\section{Uniform Estimates up to Initial Time}

In this part, we look to achieve some estimates
uniform for small time, i.e. for $t\in (0, T)$
without degeneration towards $t=0$.

The first thing comes to mind would be the
uniform $L^\infty$-estimate for the metric
potential, $u$, up to $t=0$. With the Ko\l odziej
type of estimates and even more elementary
relation between $\frac{\partial u}{\partial t}$
and $u$, naturally one wants to control the
volume form or $\frac{\partial u}{\partial t}$
up to the initial time. In order to apply the
known results from pluripotential theory to
get $L^\infty$-estimate, it is natural to require
$[\omega_0]=[\omega]$ to be semi-ample.
We focus on the case of $[\omega]$ being
big and semi-ample in this work.

Also, as in \cite{chen-tian-z}, it is reasonable
to put some restriction on the initial measure,
$\omega^n_0=(\omega+\sqrt{-1}\p\bar\p v)^
n$. If one just apply the general regularization
of the current $\omega_0$ (in \cite{blo-koj} for
example), it may not be the case that the
volume form of the smooth approximations
would have the same kind of control for the
Monge-Amp\`ere measure.

Fortunately, in sight of the discussion in
Section 1 on the uniqueness of the weak
flow, we can make proper choice of the
approximation maintaining the measure
control and also having the weak flow as
limit. More precisely, suppose the measure
$(\omega+\sqrt{-1}\partial\bar\partial v)^n$
(assumed to be $L^1$ to begin with) has
some kind of bound (for example, $upper$,
$lower$ or $L^{p>1}$-bounds), then one
can use standard process involving
partition of unity and convolution to
construct a sequence of smooth volume
forms, $\Omega_\epsilon$, having the
same kind of bound as for $(\omega+\sqrt
{-1}\partial\bar\partial v)^n$ uniformly and
converges to it as $\epsilon\to 0$ in $L^1$
space. Then we can solve the equations $
(\omega+\epsilon\omega_1+\sqrt{-1}\partial
\bar\partial v_\epsilon)^n=C_\epsilon\Omega
_\epsilon$ where $C_\epsilon$ are well
controlled and tends to $1$ after requiring
the total measure of $(\omega+\sqrt{-1}\p
\bar\p v)^n$ being positive. Taking proper
normalization for $v_\epsilon$, it would
decrease to $v$ as $\epsilon\to 0$ from
Ko\l odziej's argument as discussed in 9.6.2
of \cite{thesis}. The solutions would form a
desirable smooth approximation for the
construction of weak flow.

With all the above preparation, we look
at the situation as $t\to 0^+$ for some
cases of interests. We are still working
directly on the approximation flows, i.e.
(\ref{eq:approxi-krf}) and (\ref
{eq:approxi-skrf}), while omitting $j$
and $\epsilon_j$ for simplicity.

\subsection{Volume Upper Bound for
Semi-ample and Big $[\omega_0]$}

Let's assume the degenerate initial class
$[\omega_0]=[\omega]$ is semi-ample
and big. By the generalization of
Ko\l odziej's $L^\infty$-estimate (in \cite
{ey-gu-ze}, \cite{zzo} and \cite{demailly-pali}),
it's enough to have measure $L^{p>1}$
bound uniform up to the initial time in order
to conclude uniform $L^\infty$-bound for
$u$. In order to see this, one needs to
notice that although $\omega_t$ is changing,
if we can deal with the most degenerate one,
$\omega$, then $\omega_t$ can be treated
as $\alpha\omega+\Phi_t$ where $\alpha\in
(0,1)$ and $\Phi_t>0$ as described in \cite
{t-znote}.

\vspace{0.1in}

Now let's search for the upper bound of $\frac
{\partial u}{\partial t}$ uniform up to the initial
time under the assumption that the volume
of the initial current, $\omega_0$ has a uniform
upper bound, i.e. $L^\infty$-bound. The point
is to see whether it would be enough to imply
the necessary measure bound for small time.

The following equation comes from standard
manipulation of (\ref{eq:approxi-skrf}),
\begin{equation}
\label{eq:3}
\frac{\partial}{\partial t}\((e^t-A)\frac{\partial u}
{\partial t}-Au\)=\Delta_{\widetilde\omega_t}
\((e^t-A)\frac{\partial u}{\partial t}-Au\)+An-
\<\widetilde\omega_t, \omega-(1-A)\omega_
\infty\>
\end{equation}
for a constant $A$ to be fixed shortly. Since
the flow is driving the class into the K\"ahler
cone, it would take $A\geqslant 1$ to give
a desirable sign for the last term on the right
hand side of this equation. Unfortunately,
Maximum Principle would then give the
desirable wrong direction of control for $
\frac{\partial u}{\partial t}$ up to $t=0$.

Indeed one needs to choose some constant
$A<1$. Because $[\omega]$ is semi-ample
and big, we can choose $\omega\geqslant
0$ and a effective divisor (i.e holomorphic
line bundle) $E$ with the defining holomorphic
section $\sigma$ and a Hermitian metric $|
\cdot|$ such that for any sufficiently small
$\lambda>0$,
$$\omega_0+\lambda\sqrt{-1}\partial\bar
\partial{\rm log}|\sigma|^2>0.$$

\begin{remark}

The estimation of $\frac{\partial u}{\partial t}$
under discussion would work for $[\omega]$
being nef. and big. The only difference is one
might not have $\omega\geqslant 0$, and so
the Hermitian metric would depend on the
choice $\lambda$. We are mainly interested
in the consequential control of $u$. That's
why semi-ampleness is assumed at this
moment.

\end{remark}

Now we can reformulate (\ref{eq:3}) as follows
\begin{equation}
\begin{split}
\frac{\partial}{\partial t}\((e^t-A)\frac{\partial u}
{\partial t}-Au+\lambda{\rm log}|\sigma|^2\)
&= \Delta_{\widetilde\omega_t}\((e^t-A)\frac
{\partial u}{\partial t}-Au+\lambda{\rm log}
|\sigma|^2\)+ \\
&~~ +An-\<\widetilde\omega_t, \omega+
\lambda\sqrt{-1}\partial\bar\partial{\rm log}
|\sigma|^2-(1-A)\omega_\infty\>.
\nonumber
\end{split}
\end{equation}
By choosing $A<1$ according to the size of
$\lambda$, one can make sure the last term
on the right has a definite sign for the small
time interval.

\vspace{0.1in}

Recall that we assume $\frac{\p u}{\p t}\vline
_{\,t=0}\leqslant C$. More precisely, this means
$\omega^n_0\leqslant C\cdot\Omega$. We'll
keep using the understanding at the beginning
of this section.

Applying Maximum Principle, we have
\begin{equation}
\label{eq:deg-upper}
(e^t-A)\frac{\partial u}{\partial t}-Au+\lambda
{\rm log}|\sigma|^2\leqslant max_{X\times
\{0\}}(1-A)\frac{\partial u}{\partial t}+\lambda
\log|\sigma|^2+C\leqslant C.
\end{equation}

So one arrives at the following degenerate
upper bound
\begin{equation}
\label{eq:deg-upper2}
\frac{\partial u}{\partial t}\leqslant \frac
{-\lambda{\rm log}|\sigma|^2+C}{e^t-A}.
\end{equation}
If we can have the positive number $\frac
{\lambda}{1-A}$ small enough, then the
desired $L^p$-bound for the measure $e^
{\frac{\partial u}{\partial t}}\Omega$ can
be achieved \footnote{Remark \ref
{re:about-measure} describes why one
can ignore the term $e^u$ in the measure.}.
Notice that the constants $A$ and $\lambda$
are related. Although the constant $C$ may
also change with them, that doesn't bring any
trouble. This is indeed an assumption on
the geometry of the (effective and K\"ahler)
cones.

This degenerate upper bound (\ref{eq:deg-upper2})
gets better for large $t$. In fact, we can see in (\ref
{eq:deg-upper}) that the assumption on the initial
value can be weakened to
\begin{equation}
\label{eq:deg-upper-initial}
\frac{\partial u}{\partial t}\vline_{\,t=0}\leqslant
\frac{-\lambda\log|\sigma|^2+C}{1-A}
\end{equation}
and $\omega^n_0$ is $L^1$.

\begin{prop}

\label{prop:upper}
In the setting of Theorem \ref{th:main}, if $
\omega_0$ has $L^1$ measure satisfying
(\ref{eq:deg-upper-initial}) representing a
nef. and big class, then the weak flow
defined in Theorem \ref{th:main} would
satisfy (\ref{eq:deg-upper2}).

\end{prop}

\begin{remark}

If $\omega$ is a K\"ahler metric (as in the case
discussed in \cite{chen-tian-z}), then one doesn't
need to involve the term $\log|\sigma|^2$ in the
above estimation. The measure would have
uniform $L^\infty$-norm for short time.

\end{remark}

Although this result is not quite satisfying for the
sake of the $L^\infty$ control for $u$, it does tell
something about the convergence of the weak
flow as $t\to 0^+$ back to the initial current. We
illustrate this as follows.

The degenerate upper bound of $\frac{\p u}{\p t}$
means locally out of $\{\sigma=0\}$ \footnote{This
can be improved to be the stable base locus set
of $[\omega]$.}, $u$ is decreasing up to a term
like $Ct$ as $t\nearrow$. So by the classic result
on weak convergence (as summarized in \cite
{kojnotes}), for the weak flow,
$$\widetilde\omega^j_t\to\omega^j_0 ~\text
{weakly over}~X\setminus\{\sigma=0\}~\text
{as}~t\to 0^+, ~j=1, \cdots, n.$$

Then we can conclude the weak convergence
over $X$ for Monge-Amp\`ere measure using
the global cohomology information. Take a
sequence of strictly increasing sets exhausting
$X\setminus\{\sigma=0\}$, $\{U_k\}$ with
smooth functions $\{\rho_k\}$ supported on
$U_{k+1}$ and equal to $1$ over $U_k$. For
any non-negative smooth function $G$ over
$X$, we have
$$limit_{t\to 0^+}\int_X\rho_k\cdot G\cdot\widetilde
\omega^n_t=\int_X\rho_k\cdot G\cdot\omega^n_0.$$
Now we take supreme of the above convergence
to have
$$sup_{k}\(limit_{t\to 0^+}\int_X\rho_k\cdot G\cdot
\widetilde\omega^n_t\)=sup_{k}\(\int_X\rho_k\cdot
G\cdot\omega^n_0\)=\int_XG\omega^n_0.$$

In the mean time,
\begin{equation}
\begin{split}
\underline{limit}_{t\to 0^+}\int_XG\cdot\widetilde
\omega^n_t
&= \underline{limit}_{t\to 0^+}sup_k\(\int_X\rho_k
\cdot G\cdot\widetilde\omega^n_t\) \\
&\geqslant sup_{k}\(limit_{t\to 0^+}\int_X\rho_k
\cdot G\cdot\widetilde\omega^n_t\). \nonumber
\end{split}
\end{equation}
So we arrive at
$$\underline{limit}_{t\to 0^+}\int_XG\cdot\widetilde
\omega^n_t\geqslant \int_XG\omega^n_0.$$
The $=$ has to hold for $G\equiv 1$, and so it is
not hard to justify it for any test function $G$.
$\overline{limit}$ can be treated in the same
way. Hence we conclude
$$\widetilde\omega^n_t\to\omega^n_0 ~\text
{weakly over}~X~\text{as}~t\to 0^+.$$

\begin{corollary}

In the setting of Proposition \ref{prop:upper},
$\widetilde\omega^n_t$ converges weakly over
$X$ back to the initial current $\omega^n_0$.

\end{corollary}

\begin{remark}

The proof of this corollary can be applied
for any wedge power with $G$ replaced
by proper power of a K\"ahler metric. The
conclusion is weaker than that of \cite
{chen-tian-z} for the case there.

\end{remark}

\subsection{Volume Lower Bound}

The consideration of lower bound of volume
for small time might look a little strange, but
it gets want we want in a more elementary
way. There is also an interesting application
provided at the end.

Suppose the initial data has a positive volume
lower bound, i.e.
$$\frac{\partial u}{\partial t}\vline_{\,t=0}
\geqslant -C.$$
Together with the nef. assumption for the
main theorem, this says more or less that
$[\omega_0]$ is also big.
This lower bound can be preserved for the
approximation described at the beginning
of this section.

Recall the equation (\ref{eq:deg-upper})
appeared before
$$\frac{\partial}{\partial t}\((e^t-A)\frac
{\partial u}{\partial t}-Au\)=\Delta_
{\widetilde\omega_t}\((e^t-A)\frac{\partial
u}{\partial t}-Au\)+An-\<\widetilde\omega
_t, \omega-(1-A)\omega_\infty\>.$$
Now one choose a proper constant $A>1$
so that the last term
$$\<\widetilde\omega_t, \omega-(1-A)\omega
_\infty\>>0$$
which is possible from the cohomology picture.
Maximum Principle then gives
$$(e^t-A)\frac{\partial u}{\partial t}-Au\leqslant
max_{X\times\{0\}}(1-A)\frac{\partial u}{\partial
t}+C\leqslant C$$
in sight of the lower bound of $\frac{\p u}{\p
t}$ for the initial time. That is
$$(e^t-A)\frac{\partial u}{\partial t}\leqslant
C.$$
Hence for small time such that $e^t-A\leqslant
-C<0$, we have
$$\frac{\partial u}{\partial t}\geqslant\frac{C}
{e^t-A}\geqslant -C.$$

This automatically gives lower bound for $u$
for small time, and also the weak convergence
of $\widetilde\omega^j_t$ to $\omega^j_0$
over $X$ as $t\to 0^+$ for $j=1, \cdots, n$
from the monotonicity of $u+Ct$.

\begin{prop}

\label{prop:lower}
In the setting of Theorem \ref{th:main},
suppose $\omega^n_0$ is $L^1$ and
with a uniform positive lower bound,
then the metric potential is uniformly
bounded for small time and
$$\widetilde\omega^j_t\to\omega^j_0
~\text{weakly over}~X~\text{as}~t\to
0^+, ~j=1, \cdots, n.$$

\end{prop}

This situation occurs very naturally. For
the classic K\"ahler-Ricci flow, in finite
time singularity case, if the singular
class $[\omega_T]$ is semi-ample and
big, by parabolic Schwarz Lemma (as
in \cite{song-tian}), we have $\<\widetilde
\omega_t, \omega_T\>\leqslant C$.
So $\widetilde\omega^n_t\geqslant C
\omega^n_T$.

Here the semi-ample $[\omega_T]$
generates a map $P: X\to \mathbb
{CP}^N$ and $\omega_T=P^*\omega
_{_{FS}}$ for the standard Fubini-Study
metric $\omega_{_{FS}}$. It being big
means $P(X)$ is of the same complex
dimension as $X$.

If $F(X)$ is smooth, then the push-forward
of $\widetilde\omega_T$ would be in the
setting of Proposition \ref{prop:lower}, and
so the weak flow over $F(X)$ would weakly
converges back to the push-forward current.

{\it This simple observation makes the
picture of global weak flow on complex
surface of general type very satisfying.}

\begin{acknowledge}

The author would like to thank everyone who has supported this work
and beyond. The collaboration on a earlier work with Xiuxiong Chen
and Gang Tian has provided precious knowledge on this topic and
valuable experience. Jian Song and other people's interest and
discussion are also important. The very recent work of Jian Song and
Gang Tian (\cite{s-t-weak-flow}) considers this problem in general
algebraic geometry setting. Their results are of different flavor
from ours. Finally one can not say enough to the friendly
environment of the mathematics department at University of Michigan,
at Ann Arbor.

\end{acknowledge}

\end{document}